\documentclass[11pt,a4paper]{article}

\usepackage[cp1250]{inputenc}
\usepackage[T1]{fontenc}
\usepackage[intlimits]{amsmath}
\usepackage{amssymb}
\usepackage{amsthm}
\usepackage{geometry}
\newtheorem{tw}{Theorem}
\newtheorem{lm}{Lemma}
\newtheorem{df}{Definition}

\DeclareMathOperator*{\rot}{rot}
\DeclareMathOperator*{\di}{div}
\DeclareMathOperator*{\essup}{ess\,sup}

\DeclareMathOperator*{\tr}{Tr}

\newcommand{\nc}{\newcommand*}
\nc{\ds}{\displaystyle}
\nc{\mc}{\mathcal}
\nc{\mbb}{\mathbb}
\nc{\ee}{\end{equation}}
\nc{\be}{\begin{equation}}
\nc{\bd}{\begin{displaymath}}
\nc{\ed}{\end{displaymath}}
\nc{\barray}{\begin{array}}
\nc{\earray}{\end{array}}
\nc{\lee}{\left}
\nc{\ri}{\right}

\nc{\HH}{\mbb{H}}
\nc{\n}{\nabla}
\nc{\D}{\Delta}
\nc{\om}{\omega}
\nc{\dt}{\frac{d}{dt}}
\nc{\pdt}{\frac{1}{2}\frac{d}{dt}}
\nc{\AQ}{A\Qm}
\nc{\tf}{\tilde{F}}
\nc{\ti}{\tilde}
\nc{\h}{\hat}
\nc{\Pn}{\mathbb{P}_N}
\nc{\Qm}{Q_m}
\nc{\QQm}{Q^1_m}
\nc{\Pm}{P_m}
\nc{\PPm}{P^1 _m}
\nc{\PPn}{P^1_n}
\nc{\QQn}{Q^1_n}

\nc{\ph}{\varphi}
\nc{\aal}{\alpha}
\nc{\la}{\lambda}
\nc{\inff}{\infty}
\nc{\te}{\theta}
\nc{\bu}{\bar{u}}
\nc{\bo}{\bar{\om}}

\numberwithin{equation}{section}

\begin{document}
\setlength\arraycolsep{2pt}
\title{Finite dimensionality of 2-D micropolar fluid flow with periodic boundary conditions}
\author{P. Szopa \\
Institute of Applied Mathematics and Mechanics \\
Warsaw University, \\
Banacha 2, 02-097 Warsaw\\
szopa@duch.mimuw.edu.pl}

\maketitle

\begin{abstract}
This paper is devoted to describe the finite-dimensionality of a two-dimensional micropolar fluid flow with periodic boundary conditions. We define the notions of determining modes and nodes and estimate the number of them, we also estimate the dimension of the global attractor. Finally we compare our results with analogous results for Navier-Stokes equation.

\vspace{1cm}
\end{abstract}

\section{Introduction}
There are some heuristic as well as based on dimensional analysis arguments suggesting that the long-time behavior of a  turbulent flow is determined by a finite number of parameters.
These arguments are based on a conjecture that rapidly varying, high-wavenumber components decay as fast as they leave lower-wavenumber ingredients unaffected. It turns out from Kolmogorov's theory 
that in 3-dimensional flows only  the wavenumbers up to the cutoff value $\lambda_K = (\varepsilon/\nu^3)^{1/4}$ need to be
considered. The question is then reduced, as explained in \cite{land}, to find the number of resolution elements necessary
to be considered to describe the behavior of a fluid in a volume - say a cube of length $l_0$ on each side. The
smallest resolved distance is $l_d=1/\lambda_K$ and therefore the number of resolution elements is 
$(l_0/l_d)^3$.

A theory of Kraichnan \cite{krai}, referring to a 2-dimensional turbulent
flow, allows us to estimate the number of resolution elements to be considered as
$(l_0/\lambda_{Kr})^2$, where $\lambda_{Kr}$ is the Kraichnan length $\lambda_{Kr} = (\nu^3/\chi )^{1/6}$ and  $\chi$ is the average enstrophy dissipation rate. We refer the reader to \cite{doering} for more detailed discussion of turbulent length scales. 

The notion of determining modes arise naturally during consideration of the Fourier decomposition of a flow. There are some results concerning determining modes in the context of 2-dimensional Navier-Stokes equations. Foias and Prodi showed in \cite{foip} that if a number of Fourier modes of two different solutions have the same asymptotic behavior, then whole solutions have the same asymptotic behavior too. Subsequent works have been aimed at estimating how many low modes are necessary to determine the behavior of a flow. The most recent results are in \cite{foi2} for the case of no-slip boundary conditions and in \cite{jones} for the case of periodic boundary conditions. The number of determining modes for 2-D micropolar fluid flow with no-slip boundary conditions was estimated in \cite{szopad}.

In many practical situations, for instance physical experiments, data are collected from measurements in finite points in a domain of a flow. Therefore occurs a question how many measurement points are necessary to determine the long-term behavior of a flow. This leads to the notion of determining nodes, which was introduced by Foias and Temam in \cite{foitem}. The most recent estimate for the lowest number of determining nodes for Navier-Stokes equation in space-periodic case was derived in \cite{jones}.

Another approach to describe the asymptotic behavior of a flow with a finite number of parameters is to use a global attractor.
For every trajectory $u$ in the phase space we can choose, due to "Shadowing Lemma" (c.f. \cite{rob}), a trajectory $u_A$ lying on the attractor, that is arbitrary close in an interval of time, that is $|u(t) -u_A(t)| < \varepsilon$ for $t\in(t_0,t_1)$. On the other hand the global attractor has finite Hausdorff and fractal dimension, so we can parametrize it by a finite number of parameters (c.f. \cite{mane}, \cite{hunt}, \cite{foio}). Therefore we can describe approximately the long-term behavior of a flow by a finite number of parameters.

There are many results concerning the dimension of attractor for Navier-Stokes equation with a variety of boundary conditions cf. ie.  periodic boundary conditions in \cite{tem2}, pipe-like domain with arbitrary inflow at infinity in \cite{muchasad}.
The dimension of attractor for micropolar fluid equations with various boundary conditions was estimated in \cite{luk1}, \cite{bouluk}.

In this paper we will consider all of mentioned above ways of determining the long-time behavior of a micropolar fluid flow by finite number of parameters. We estimate the lowest number of determining modes and nodes, and the dimension of the global attractor.

We will consider the micropolar fluid equations, which in two-dimensional case have form (c.f. \cite{luk1}) 
\be \label{eq1}
\frac{\partial u}{\partial t} - (\nu + \nu_{r}) \D u + (u \cdot\n) u + \n p = 2 \nu_{r} \rot \om + f,
\ee 
\be\label{eq2} 
\di u =0, 
\ee 
\be\label{eq3} 
\frac{\partial \om }{\partial t} - \alpha \D\om  + (u \cdot \n) \om +4\nu_r \om  = 2 \nu_r \rot u +g, 
\ee
where $u=(u_1,u_2)$ is the velocity field, $p$ is the pressure and $\om$ is the microrotation field interpreted as the angular velocity of particles. In 2-dimensional case we assume that the axe of rotation of particles is perpendicular to $x_1,x_2$  plane. The fields $f=(f_1,f_2)$ and $g$ are the external forces and moments respectively. Positive constants $\nu,\nu_r, \aal$ are the viscosity coefficients and

\bd 
\rot u = \frac{\partial u_2}{\partial x_1} - \frac{\partial u_1}{\partial x_2}, \qquad 
\di u = \frac{\partial u_1}{\partial x_1} + \frac{\partial u_2}{\partial x_2}, \qquad
\rot \om = \left(\frac{\partial \om}{\partial x_2}, -\frac{\partial \om}{\partial x_1} \right).
\ed 
We supplement these equations with initial condition 
\bd u(x,0)=u_0(x), \quad \om(x,0)=\om_0(x) 
\ed 
and periodic boundary conditions 
\bd 
u(x+Le_i,t)=u(x,t),\quad \om(x+Le_i,t)=\om(x,t) 
\quad \forall x\in\mbb{R}^2 \quad \forall t>0, 
\ed 
where $e_1,e_2$ is Euclidean base of $\mbb{R}^2$ and $L$ is the period in the i-th direction; $Q=(0,L)^2$ is the square of the period.

The existence and uniqueness of solutions for this model as well as existence of the global attractor was proved in \cite{szopae}. We also assume (like in \cite{szopae}) that the space averages of $u,\om,f$ and $g$ vanish.

This paper is organised as follows: 

In section \ref{spaces} we introduce function spaces and operators used throughout this article.

Section \ref{ape} presents some \emph{a priori} estimates we need in the sequel. 

In section \ref{modes} we estimate the number of determining modes. 

The number of determining nodes is estimatef in section \ref{nodes}.

Section \ref{dimension} is devoted to recall the notions of fractal and Hausdorff dimension and to estimate the dimension of global atractor.

\section{Mathematical setting of the problem}\label{spaces}
In this section we introduce function spaces and trilinear forms $b$ and $b_1$, Stokes and $-\D$ operators and operator $\rot$ which are necessary in our consideration. 
\\ \\
\noindent \emph{Function spaces.}\\
We will denote by $\mbb{X}$ the space $X \times X$ with standard product norm .

$L^q$ is the usual Lebesque's space $L^q(Q)$ for $q\in [1,\infty]$. We denote the scalar product in $L^2$ by $(\cdot, \cdot)$ and the norm in $L^2$ by $|\cdot|$ when it doesn't produce confusion.

$H^m$ is usual Sobolev spaces $H^m(Q)$ for $m\in \mbb{N}$.

By $H^m _p(Q), \, m\in \mbb{N}$ we denote the space of real functions which are in $H^m _{loc}(\mbb{R}^n)$ and periodic with the period $Q$:
$ u(x+Le_i) = u(x) \quad i=1,2. $

For an arbitrary $m \in \mbb{N}, \, H^m _p(Q)$ is a Hilbert space with the scalar product
\bd
(u,v)_m = \sum _{|\alpha| \leq m} \int_Q D^\alpha u(x) D^\alpha v(x)\, dx 
\ed 
and the norm induced by it. The functions in $H^m _p(Q)$ are explicitly characterized by their Fourier series
expansion 
\bd H^m _p (Q)= \left\{u \colon u=\sum_{k \in \mbb{Z}^n} u_k e^{2 i \pi k/L \cdot x},\, \bar {u}_k = u_{-k},\,
         |u|_m =\sum_{k \in \mbb{Z}^n} |k|^{2m}|u_k|^2 <\infty \right\},
\ed
where $k/L = (k_1/L,k_2/L)$, the norm $|u|_m$ is equivalent to the norm 
$\{\sum_{k\in\mbb{Z}^n} (1+|k|^{2m})|u_k|^2\}^{1/2}$. We also set
\bd
\dot {H}^m _p(Q)=\{u \in H^m_p(Q) ,c_0=0\}.
\ed

$H$ and $V$ are the divergence-free subsets  of $\dot{\HH}^0_p(Q)$ and $\dot{\HH}^1_p(Q)$, respectively. We equip $V$ with the scalar product and the Hilbert norm
\bd
 ((u,v))=\sum_{i=1} ^n \left(\frac{\partial u}{\partial x_i},\frac{\partial v}{\partial x_i} \right),
        \qquad ||u||=\{((u,u))\}^{1/2}.
\ed
This norm is equivalent to the norm induced by $H^1_p(Q)$, and $V$ is a Hilbert space for this scalar product.

One can check that $\dot H^{-m}_p$ is the dual space to $\dot H^m_p$, we also denote the dual space to $V$ as $V'$.

Let $\mc{H}$ and $\mc{V}$ denote $H \times \dot H_p ^0$ and $V \times \dot H_p ^1$, respectively, with standard product norms.

$L^q(0,T;X)$ is the space of strongly measurable functions $u\colon(0,T) \to X$, where $X$ is a Banach space, with the followning norm
\bd
||u||_{L^q(0,T;X)}=\left\{\begin{array}{ll}
                    \ds{\left(\int_0 ^T ||u(t)||_X ^q \, dt\right)^{1/q},} & \quad1\leq q < \infty, \\
                    \ds{\essup_{t\in(0,T)} ||u(t)||_X,} & \quad q=\infty . \end{array}   \right.
\ed

C([0,T],X) is the space of continuous functions $u\colon(0,T) \to X$, where $X$ is a Banach space, with the usual norm.
\\
\\
\noindent \emph{Stokes and $-\D$ operators.}

Let us consider Stokes problem, which is obtained from the Navier-Stokes equation by neglecting all time-dependent and nonlinear terms. 
Thus it is the following one: for given $f\in \dot{\HH}_p^0$ or $\dot{\HH}^{-1}_p$, find $u \in \dot{\HH}_p^1$ and $p\in L^2$ such that
\be \label{stokes1}
-\D u + \n p=f, \qquad \n \cdot u =0.
\ee
It is shown (cf., e.g.  \cite{rob,tem}) that the Stokes operator $A$ associated with space periodicity condition is simply $-\D$ (provided that $f \in H$) with domain $D(A)=\dot H_p^2 \cap H$.  $A$ is one-to-one map form $D(A)$ onto $H$.

The operator $A^{-1}$ is linear, continuous from $H$ into $D(A)$. Since the injection of $D(A)$ in $H$ is compact, we can consider $A^{-1}$ as a compact operator in $H$. $A^{-1}$ is also self-adjoint as an operator in $H$. Hence it possesses a sequence of eigenfunctions $w_j$, $j\in\mbb{N}$ which form an orthonormal basis of $H$
\bd
\barray{l}
Aw_j=\la_j w_j, \qquad w_j \in D(A),\\
0<\la_1 \leq \la_2 \leq \ldots, \qquad \la_j \to \infty  \quad \textrm{for} \quad j\to \infty.
\earray
\ed
The  operator  $-\D$ has the same properties as the Stokes operator. The eigenvalues are the same but the eigenfunctions are different, because we consider $\om$ to be in $\mbb{R}$, and we denote them by $\rho_k$. For notation consistence we write $A_1$ instead of $-\D$.
 
We can express every element $u\in H$ and $\om \in \dot H^0_p$ as
\begin{displaymath}
u(x,t)= \sum _{k=1} ^{\infty} u _k (t) w_k(x), \qquad
\om(x,t)= \sum _{k=1} ^{\infty} \om _k (t) \rho_k(x).
\end{displaymath}
The Galerkin projectors corresponding to the first $m$ modes are:
\begin{displaymath}
\Pm u(x,t)= \sum _{k=1} ^{m} u_k (t) w_k(x) \qquad
\PPm \om_i(x,t)= \sum _{k=1} ^{m} \om_k (t) \rho_k(x).
\end{displaymath}
 We also denote projections onto modes higher than $m$ by $\Qm$ and $\QQm$ 
\bd
\Qm u(x,t)= \sum _{k=m+1} ^{\infty} uk (t) w_k(x) \quad
\QQm \om(x,t)= \sum _{k=m+1} ^{\infty} \om_k (t) \rho_k(x).
\ed

\noindent \emph{Trilinear forms.}\\
We define the trilinear forms $b$ and $b_1$ as follows
\bd
b(u,v,w) = \sum_{i,j=1} ^2 \int_Q u_i \frac{\partial v_j}{\partial x_i} w_j \,dx
\ed
for all $u,v,w \in V$ and
\bd
b_1(u,\om,\psi) = \sum_{i=1} ^2 \int_Q u_i \frac{\partial \om}{\partial x_i} \psi \,dx
\ed
for all $u \in V,$ and all scalar functions $\om, \psi \in \dot{H}^1_p(Q)$.
The forms $b$ and $b_1$ have the orthogonality property
\bd
b(u,v,v)=0,\quad b_1(u,\om,\om)=0.
\ed
In the 2-dimensional space-periodic case the form $b$ has one more orthogonality property \cite{tem}
\be \label{ortho}
b(u,u,Au)=0 \quad \forall u\in D(A),
\ee
which the form $b_1$ does not have -- it is not true that $b_1(u,\om,A_1\om)= 0$ for all $u \in D(A), \om \in D(A_1)$.
The lack of this orthogonality property causes that the a priori estimates, we obtain in the sequel, are more involved than analogous estimates for NSE with periodic boundary conditions.

We do some estimates of the forms $b$ and $b_1$ using the Ladyzhenskaya inequality \cite{lad}
\bd
||u||_{L^4} \leq \lee(\frac{6}{\pi}\ri)^{1/4}|u|^{1/2}||u||^{1/2}, \quad \textrm{for all } u\in \dot H_p^1,
\ed
and the Agmon inequality \cite{agmon, ilyin}
\bd
||u||_{L^\infty} \leq \frac{1}{\sqrt{\pi}}|u|^{1/2}|Au|^{1/2}, \quad \textrm{for all } u\in D(A).
\ed
The most relevant estimate of the constants are obtained in \cite{ilyin}. We also use the H\"older inequality:
\be \label{b4}
\begin{array}{ll}
|b(u,v,w)|      \leq c_1 |u|^{1/2}\|u\|^{1/2}\|v\| \cdot |w|^{1/2} ||w||^{1/2}   & u,v,w \in V,\\
|b(u,v,w)|      \leq c_1 |u|^{1/2}||u||^{1/2}||v||^{1/2}|Av|^{1/2}|w|  & u\in V, \, v,w \in D(A),\\
|b(u,v,w)|      \leq c_1 |u|^{1/2}|Au|^{1/2}||v||\cdot|w|  & u\in D(A),\,v\in V, \, w \in H,\\
|b_1(u,\om,\psi)| \leq c_1 |u|^{1/2}||u||^{1/2}|\psi|^{1/2}||\psi||^{1/2}||\om|| &  u,\om,\psi \in \dot H^1 _p, \\
|b_1(u,\om,A_1\psi)|\leq c_1 |u|^{1/2}|Au|^{1/2}||\om||\cdot|A_1\psi| &  u\in D(A),\, \om\in \dot{H}_p^1, \,
                                                                 \psi \in D(A_1),\\
|b_1(u,\om,A_1\psi)|\leq c_1 |u|^{1/2}||u||^{1/2}||\om||^{1/2}|A_1\om|^{1/2}|A_1\psi| &  u \in V,\,\om,\psi \in D(A_1)
\end{array}
\ee
for an appropriate constant $c_1$.

The operator $\rot$ has the following properties,
\begin{eqnarray}
 \ds{ \int_Q \rot u \cdot \om \, dx}&=&\ds{ \int_Q \rot \om \cdot u \, dx, } \label{id1} \\[1.5ex]
 \ds{\int_Q |\rot \om|^2 \,dx }&=&\ds{ \int_Q |\n \om|^2 \, dx, }\label{id2} \\[1.5ex]
 \ds{\int_Q |\rot u|^2 \,dx}& =&\ds{ \int_Q |\n u|^2 \, dx, }\label{id3}
\end{eqnarray}
for all $u \in V$ and $\om \in \dot H^1 _p$. 

\section{\emph{A priori} estimates} \label{ape}
In this section we derive some \emph{a priori} estimates, which we use in the sequel. We will obtain these in the terms of asymptotic strength of forces and moments. To this end we set
\be\label{F}
\tf = \limsup_{t\to\infty} \left(|f(t)|^2 + |g(t)|^2\right)^{1/2}
\ee
and
\be \label{k}
k_1 = \min\{\nu,\alpha\},\quad k_2=k_1\la_1.
\ee
All inequalities, apart form the latter, we use throughout this section come from \cite{szopae}.
Let us consider the following one
\be\label{e1}
\ds{\frac{d}{dt} (|u(t)|^2 + |\om(t)|^2) }+ k_2 \left(|u(t)|^2 + |\om(t)|^2\right)
                \leq k_2 ^{-1} \left(|f(t)|^2 +|g(t)|^2\right).
\ee 
Integrating it with respect to $t$ in the interval $(0,t)$ we obtain in particular 
\be\label{e2}
\begin{array}{rl}
|u(t)|^2 +& |\om(t)|^2  \leq \ds{|u_0|^2 +|\om_0|^2 + k_2 ^{-1} \int_0^t \left(|f(s)|^2 +|g(s)|^2\right)\, ds  }\\
                \leq& |u_0|^2 +|\om_0|^2
                     \ds{+ k_2 ^{-1} \left(||f||_{L^2(0,T;H)} ^2 +||g||_{L^2(0,T;\dot H_p ^0)}^2\right)}\\
                \leq&  |u_0|^2 +|\om_0|^2 + \ds{k_2 ^{-1} \left(||f||_{L^2(0,\infty;H)} ^2
                        +||g||_{L^2(0,\infty;\dot H_p ^0)}^2\right)},
\end{array}
\ee
which implies the uniform bound of the norm of solution in $\mc{H}$.

We apply Gronwall's inequality to \eqref{e1} to obtain an uniform with respect to initial conditions bound of the norm of solutions for large times. We proceed as follows
\begin{align*}
|u(t)|^2 +|\om(t)|^2    \leq& e^{-k_2(t-t_0)}\Bigg\{|u(t_0)|^2+|\om(t_0)|^2\\
                        &+k_2^{-1}\int_{t_0}^t \left(|f(s)|^2 +|g(s)|^2\right)e^{k_2(s-t_0)} \, ds\Bigg\} \\
                \leq& e^{-k_2(t-t_0)}\Bigg\{|u(t_0)|^2+|\om(t_0)|^2\\
                        &+k_2^{-1}\int_{t_0}^t \left(||f||^2_{L^\infty(t_0,t;H)}
                        +||g||^2_{L^\infty(t_0,t;\dot H_p ^0)}\right)e^{k_2(s-t_0)}\, ds\Bigg\} \\
                \leq& e^{-k_2(t-t_0)}\Bigg\{|u(t_0)|^2+|\om(t_0)|^2\\
                        &+k_2^{-2}\left(||f||^2_{L^\infty(t_0,t;H)}+||g||^2_{L^\infty(t_0,t;\dot H_p^0)}\right)
                        \left(e^{k_2(t-t_0)}-1\right) \Bigg\} \\
                \leq& e^{-k_2(t-t_0)}\left(|u(t_0)|^2+|\om(t_0)|^2\right)\\
                        &+k_2^{-2} \left(1-e^{-k_2(t-t_0)}\right)
                        \left(||f||^2_{L^\infty(t_0,t;H)}+||g||_{L^\infty(t_0,t;\dot H_p ^0)}^2\right), \\
\end{align*}
hence for $t_0$ and $t$ large enough
\be\label{e4}
|u(t)|^2 +|\om(t)|^2 \leq \frac{2}{k_2^2}\tf^2.
\ee
To estimate the average of the square of the norm of solutions in $\mc{V}$ we use the inequality 
\be\label{e5}
\frac{d}{dt} (|u(t)|^2 + |\om(t)|^2) + k_1 \left(||u(t)||^2 + ||\om(t)||^2\right)
                \leq k_2 ^{-1} \left(|f(t)|^2 +|g(t)|^2\right).
\ee
By integration we obtain
\begin{align*}
|u(t+T)|^2 +|\om(t+T)|^2 &+ k_1\int_t ^{t+T}\lee(||u(s)||^2 +||\om(s)||^2\ri)\, ds \\
                &\leq k_2^{-1} \int_t ^{t+T}\lee(|f(s)|^2+|g(s)|^2\ri)\, ds + |u(t)|^2 +|\om(t)|^2.
\end{align*}
Because $|u(t)|^2 +|\om(t)|^2$ is uniformly bounded with respect to $t$ (c.f. \eqref{e2}), then for $t$ and $T$ large enough we have
\be\label{e6}
\begin{array}{rl}
\ds{\frac{1}{T}\int_t ^{t+T}\lee||u(s)||^2 +||\om(s)||^2\ri) \,ds} \leq& \ds{(k_1k_2)^{-1}\frac{1}{T}\int_t ^{t+T}|f(s)|^2+|g(s)|^2\, ds 
                                        + \frac{1}{T}\left(|u(t)|^2 +|\om(t)|^2\right) }\\[1.5ex]
                                    \leq&\ds{ \frac{2}{k_1k_2} \tf^2.}
\end{array}
\ee

In order to derive two more estimates we consider the inequality 
\be\label{e7}
\begin{array}{rl}
\ds{\dt\big(||u(t)||^2 +}&\ds{||\om(t)||^2\big) + \frac{k_1}{2}\lee(|Au(t)|^2+|A_1\om(t)|^2\ri)} \\[1.5ex]
    \leq& \ds{ \left(\frac{2C}{\alpha^2 \nu}|u(t)|^2||\om(t)||^2
                    +\frac{8\nu_r}{\alpha}\right)\lee(||u(t)||^2+||\om(t)||^2\ri)}  \\[1.5ex]
    &\ds{+\frac{2}{k_1}\left(|f(t)|^2+|g(t)|^2\right)}.
\end{array}
\ee
Let us set
\begin{align*}
y(t)        &=||u(t)||^2+||\om(t)||^2,\\
\ti{g}(t)   &=\left(\frac{2C}{\alpha^2 \nu}|u(t)|^2||\om(t)||^2+\frac{8\nu_r}{\alpha}\right),\\
\ti{h}(t)   &=\frac{2}{k_1}\left(|f(t)|^2+|g(t)|^2\right),
\end{align*}
then we can rewrite \eqref{e7} in the form
\bd
\frac{dy}{dt} \leq \ti{g}y +\ti{h}.
\ed
We check assumptions of uniform Gronwall's lemma (c.f. \cite{tem2}). For $t$ and $r$ large enough we have
\begin{align*}
\int_t^{t+r}\ti{g}(s) \,ds  \leq&  \int_t^{t+r} \lee\{\frac{2C}{\alpha^2\nu}\lee(|u(s)|^2+|\om(s)|^2\ri)
                    \lee(||u(s)||^2+||\om(s)||^2\ri) +\frac{8\nu_r^2}{\alpha} \ri\}\,ds \\
                \leq& \frac{8\nu_r^2 r}{\aal} + \frac{4C}{\aal^2 \nu k_2^2} \tf^2
                    \int_t^{t+r} \lee(||u(s)||^2+||\om(s)||^2\ri)\, ds \\
                \leq& \frac{8\nu_r^2 r}{\aal} + \frac{8Cr}{\aal^2\nu k_1 k_2^3} \tf^4 \equiv a_1,\\
\int_t^{t+r}\ti{h}(s) \,ds  \leq& \frac{3r}{k_1} \tf^2 \equiv a_2,\\
\int_t^{t+r}y(s) \,ds   \leq& \frac{2r}{k_2^2} \tf^2 \equiv a_3.
\end{align*}
Therefore by uniform Gronwall's lemma we obtain 
\be\label{e7a}
||u(t)||^2+||\om(t)||^2 \leq \frac{2 + 3k_2 r}{k_1k_2}\tf^2
            \exp\lee(\frac{8\nu_r^2 r}{\aal} + \frac{8C r}{\aal^2 \nu k_1 k_2^3} \tf^4\ri),
\ee
for all $t>t_0+r$, $t_0$ such that above estimates holds.
Let us denote for notation simplicity
\bd
\h{c}_1=\frac{2 + 3k_2r}{k_1k_2} , \quad \h{c}_2=\frac{8\nu_r^2r}{\aal},
        \quad \h{c}_3=\frac{8Cr}{\aal^2 \nu k_1 k_2^3},
\ed
then \eqref{e7a} becomes
\be \label{e8}
||u(t)||^2 + ||\om(t)||^2 \leq \h{c}_1 \tf^2 \exp\lee(\h{c}_2 + \h{c}_3 \tf^4 \ri).
\ee

Now we want to derive the estimate on the average of square of the norm of solution in $D(A) \times D(A_1)$. Integrating \eqref{e7} on the interval $(t,t+T)$ with respect to $t$ we get
\begin{align*}
||u(t+T)||^2 +&||\om(t+T)||^2 - ||u(t)||^2 - ||\om(t)||^2 \\
            &+ \frac{k_1}{2}\int_t^{t+T}\lee(|Au(s)|^2 +|A_1\om(s)|^2 \ri) \, ds \\
        \leq& \int_t^{t+T}\Bigg\{ \lee(\frac{2C}{\alpha^2 \nu}|u(t)|^2||\om(t)||^2
                    +\frac{8\nu_r}{\alpha}\ri)\lee(||u(t)||^2+||\om(t)||^2\ri)\\
    &+\frac{2}{k_1}\left(|f(t)|^2+|g(t)|^2\right)\, \Bigg\}ds,
\end{align*}
hence
\begin{align*}
\frac{1}{T}\int_t^{t+T} \big(|Au(s)|^2 +&|A_1\om(s)|^2\big) \,ds \leq  \frac{1}{T} \frac{2}{k_1}\lee(||u(t)||^2 + ||\om(t)||^2\ri) 
            + \frac{4}{k_1^2} \frac{1}{T}\int_t^{t+T}\lee(|f(s)|^2+|g(s)|^2\ri) \, ds \\
            &+ \frac{2}{k_1} \frac{1}{T} \int_t^{t+T} \Bigg(\frac{2C}{\alpha^2\nu}\lee(|u(s)|^2+|\om(s)|^2\ri)
            \cdot  \lee(||u(s)||^2+||\om(s)||^2\ri) +\frac{8\nu_r^2}{\alpha} \Bigg) \\
             & \cdot       \lee(||u(s)||^2+||\om(s)||^2\ri) \,ds.
\end{align*}
Because solutions are uniformly bounded in the $\mc{V}$ norm for large $t$, then for $t$ and $T$ large enough
\be \label{Au-norm}
\barray{rl}
\ds{\frac{1}{T}\int_t^{t+T} \big(|Au(s)|^2 +}&\ds{|A_1\om(s)|^2\big) \, ds \leq  \frac{5}{k_1^2}\tf^2
            + \frac{16\nu_r^2}{\aal k_1}\frac{1}{T}\int_t^{t+T}\lee(||u(s)||^2+||\om(s)||^2\ri) \, ds }\\
        &\ds{+\frac{4C}{\aal^2 \nu k_1} \frac{1}{T}
            \int_t^{t+T}\lee(|u(s)|^2+|\om(s)|^2\ri)\lee(||u(s)||^2+||\om(s)||^2\ri)^2 \, ds  }\\
        \leq &\ds{ \lee(\frac{5}{k_1^2} + \frac{32\nu_r^2}{\aal k_1^2 k_2}\ri) \tf^2 
            +\frac{8C}{\aal^2 \nu k_1 k_2^2} \tf^2 }\\
			&\ds{\cdot \frac{1}{T} \int_t^{t+T}\lee(||u(s)||^2+||\om(s)||^2\ri) \lee(||u(s)||^2+||\om(s)||^2\ri) \, ds}\\[1.5ex]
        \leq & \ds{\lee(\frac{5}{k_1^2} + \frac{32\nu_r^2}{\aal k_1^2 k_2}\ri) \tf^2 +
                \frac{16C\h{c}_1}{\aal^2\nu k_1^2 k_2^3} \tf^6 \exp \lee(\h{c}_2 + \h{c}_3\tf^4\ri)}.
\end{array}
\ee
We need one more to estimate the number  of determining modes in terms of $H^{-1}$ norm of forces and moments.
Taking the scalar product of \eqref{eq1} with $u$ in $H$ we obtain
\be \label{est10}
\frac{1}{2} \frac{d}{dt} |u|^2 + (\nu+\nu_r)||u||^2 = 2\nu_r(\rot \om,u) +(f,u)
\ee
because $b(u,u,u)=0$. We estimate the terms of the RHS of \eqref{est10} as follows
\be \label{est20}
\begin{array}{rl}
2\nu_r(\rot \om,u) &\ds{= 2\nu_r(\om,\rot u) \leq 2\nu_r|\om| \cdot ||u|| \leq 2\nu_r |\om|^2 + \frac{\nu_r}{2} ||u||^2}, \\
(f,u) &\ds{\leq ||f||_{\mbb{H}^{-1}}||u||_H \leq \frac{\nu}{2}||u||^2 + \frac{1}{2\nu}  ||f||_{\mbb{H}^{-1}}^2 }.
\end{array}
\ee
We treat \eqref{eq3} in an analogous way. We multiply it by $\om$, integrate over $Q$ and we have
\be \label{est30}
\frac{1}{2} \frac{d}{dt} |\om|^2 + \aal ||\om||^2 + 4\nu_r|\om|^2 =2\nu_r(\rot u,\om) + (g,\om).
\ee
The terms of the RHS of \eqref{est30} are estimaded
\be\label{est40}
\begin{array}{rl}
2\nu_r(\rot u,\om) & \ds{\leq 2\nu_r|\om|^2 + \frac{\nu_r}{2} ||u||^2, }\\[1.5ex]
(g,\om) \leq ||g||_{H^{-1}} ||\om|| & \ds{\leq \frac{\aal}{2}||\om||^2 + \frac{1}{2\aal} ||g||_{H^{-1}}^2. }
\end{array}
\ee
Adding equations \eqref{est10} and \eqref{est30}, using estimates \eqref{est20} and \eqref{est40} we arrive at
\be\label{est50}
\frac{d}{dt}\lee(|u|^2 +|\om|^2 \ri) +k_1\lee(||u||^2 +||\om|^2 \ri) \leq \frac{1}{k_1} \lee(||f||_{\mbb{H}^{-1}}^2 + ||g||_{H^{-1}}^2\ri).
\ee
Let us notice that \eqref{est50} looks similar to \eqref{e5}. Therefore, proceeding in the same way we obtain
\be\label{est60}
\frac{1}{T} \int_t ^{t+T} (||u(s)||^2 + ||\om(s))||^2 )\, ds \leq \frac{2}{k_1^2} \tf_{-1}^2
\ee
for $t$ and $T$ large enough.
Finally we have all \emph{a priori} estimates we need in the sequel.

\section{Determining modes}\label{modes}
In this section we introduce the notion of determining modes for micropolar fluids and estimate the number of them in terms of parameters of the model. We follow the method described in \cite{foi}.

Let us consider two solutions of 2-D micropolar fluid equations $(u_1,\om_1)$ and $(u_2,\om_2)$
($u_i=(u_i ^1(x,t), u_i ^2(x,t))$ and $\om_i = \om_i(x,t)$) corresponding to two possibly different pairs of external forces and moments $(f_1,g_1)$ and $(f_2,g_2)$ respectively. More explicitly, $(u_1,\om_1)$ and $(u_2,\om_2)$ satisfy the equations
\be \label{r11}
\frac{\partial u_{1}}{\partial t} - (\nu + \nu_{r}) \D u_1 + (u_1 \cdot \n) u_1 + \n p = 2 \nu_{r} \rot \om_1 + f_1,
\ee
\be\label{r12}
\di u_1 = 0,
\ee
\be\label{r13}
\frac{\partial \om _1}{\partial t} - \alpha \D \om_1  + (u_1 \cdot \n) \om_1 +4\nu_r \om_1  = 2 \nu_r \rot u_1 +g_1
\ee
and
\be \label{r21}
\frac{\partial u_{2}}{\partial t} - (\nu + \nu_{r}) \D u_2 + (u_2 \cdot \n) u_2 + \n q = 2 \nu_{r} \rot \om_2 + f_2,
\ee
\be\label{r22}
\di u_2 = 0,
\ee
\be\label{r23}
\frac{\partial \om _2}{\partial t} - \alpha\D \om_2  + (u_2 \cdot \n) \om_2 + 4\nu_r \om_2  = 2 \nu_r \rot u_2+g_2,
\ee
with corresponding pressure terms $p=p(x,t)$ and $q=q(x,t)$. The boundary conditions are periodic for both problems.
It is assumed, that external forces $f_1$, $f_2$ and moments $g_1$, $g_2$ have the same asymptotic behavior for large time, that is,
\be
\int _Q \lee( |f_1(x,t) - f_2(x,t)|)_{\mbb{H}^{-1}}^2 + |g_1(x,t) - g_2(x,t)|_{H^{-1}} ^2 \ri)\, dx \to 0 
		\textrm{ as } t \to \infty. \label{s1}
\ee

\begin{df}\label{mfmod}
The first $m$ modes associated with $\Pm$ and $\PPm$ are called \emph{determining modes} if the condition
\be
\int _Q \big(|\Pm u_1(x,t) - \Pm u_2(x,t)|^2 +|\PPm\om_1(x,t) - \PPm \om_2(x,t)|^2 \big)\, dx \to 0  \textrm{ as } t \to \infty \label{md1}
\ee
and the condition for forces and moments \eqref{s1} imply
\be \label{md3}
\int _Q \big(|\Qm u_1(x,t)-\Qm u_2(x,t)|^2+|\QQm\om_1(x,t)-\QQm\om_2(x,t)|^2\big) \, dx \to 0 \textrm{ as } t \to \infty.
\ee
\end{df}

We will give the estimate of the value of the number $m$ in terms of asymptotic strength of the forces and moments measured in terms of their $L^2$ and $H^{-1}$ norm, that is 
\bd 
\barray{l}
\ds{\tf= \limsup _{t \to \infty} \left( |f_1(t)|^2 + |g_1(t)|^2 \right)^{1/2}, }\\ 
\ds{\tf_{-1} = \limsup_{t\to\infty} \lee( ||f_1(t)||_{\mbb{H}^{-1}}^2 + ||g_1(t)||_{H^{-1}}^2 \ri)^{1/2}.} 
\earray
\ed

\begin{tw}
Let $f_i\in L^\infty(0,T;\mc{H}),\, g_i \in L^\infty(0,T;L^2)$ for $i=1,2$. If the forces and moments satisfy condition \eqref{s1}, then the first $m$ modes are determining in the sense of Definition \ref{mfmod} provided that
\bd 
m \geq \frac{16 \nu_r^2}{d\la_1\aal k_1}  + \frac{8c_1^2}{d\la_1k_3k_1^3} \tf^{2}_{-1}.
\ed
\end{tw}
\begin{proof}
This estimate is similar to that obtained in \cite{szopad} but in this paper we have relaxed the convergence of forces and moments to be only in $H^{-1}$. Moreover, the obtained estimate is in terms of their $H^{-1}$ norm. The same argumentation works in the case of no-slip boundary coditions.

Our estimate is based on the following generalization of the classical Gronwall's lemma (c.f. \cite{foi}).
\begin{lm} \label{gron}
Let $\gamma =\gamma(t)$ and $\beta = \beta(t)$ be locally integrable real-valued functions on $[t_0, \infty)$ that satisfy the following conditions for some $T >0$:
\be \label{l1}
\liminf _{t \to \infty} \frac{1}{T} \int _t ^{t+T} \gamma(\tau) \, d\tau >0,
\ee
\be \label{l2}
\limsup _{t \to \infty} \frac{1}{T} \int _t ^{t+T} \gamma^{-}(\tau) \, d\tau < \infty,
\ee
\be \label{l3}
\lim _{t \to \infty} \frac{1}{T} \int _t ^{t+T} \beta^{+}(\tau) \, d\tau =0,
\ee
where $\gamma^-(t) = \max\{-\gamma(t),0\}$ and $\beta ^+ (t) = \max \{\beta(t),0\}$. Suppose that $\xi = \xi(t)$ is an absolutely continuous nonnegative function on $[t_0,\infty)$ that satisfies the following inequality almost everywhere on $[t_0,\infty)$:
\bd
\frac{d\xi}{dt} + \gamma\xi \leq \beta.
\ed
Then $\xi(t) \to 0$ as $t \to \infty$.
\end{lm}

Writing the equations of micropolar fluid in functional form for a pair of solutions $(u_1,\om_1)$ and $(u_2, \om_2)$ and subtracting them we find
\begin{eqnarray}
u_t + (\nu + \nu_r)A u + B(u,u_1) + B(u_2,u)  &=& 2\nu_r \rot \om +f,    \label{dmh10} \\
\om_t + \alpha A_1 \om +B_1(u,\om_1)+ B_1(u_2,\om)  + 4\nu_r \om &=& 2\nu_r \rot u +g, \label{dmh20}
\end{eqnarray}
where $u=u_1-u_2$, $\om=\om_1-\om_2$, $f=f_1-f_2$ and $g=g_1-g_2$.
First we will deal with the equation \eqref{dmh10}, multiplying it by $\Qm u$ and integrating over $Q$ we obtain
\be \label{dmh30}
\begin{array}{rl}
\ds{\frac{1}{2}\frac{d}{dt}|\Qm u|^2 + (\nu +\nu_r)\|\Qm u \|^2 + b(u,u_1,\Qm u) }&+ b(u_2,u,\Qm u) \\
							&= 2\nu_r (\rot \om, \Qm u) + (f,\Qm u).\\
\end{array}
\ee
We estimate the linear terms of the RHS of \eqref{dmh30} as follows
\bd
\begin{array}{rl}
(f,\Qm u) &\leq ||f||_{\mbb{H}^{-1}} ||\Qm u ||_{\mbb{H}^1} \leq \ds{\frac{\nu+2\nu_r}{4} ||\Qm u||^2 + \frac{1}{\nu+2\nu_r}||f||^2 _{\mbb{H}^{-1}}  },
																																																								\\[1.5ex]
2\nu_r(\rot \om, \Qm u) &= 2\nu_r[(\Pm \rot \om, \Qm u) + (\Qm \rot \om,\Qm u)] \leq 2\nu_r |\Qm \rot \om| \cdot |\Qm u| \\[1.5ex]
										& \ds{\leq \frac{\aal}{8} ||\QQm \om||^2 + \frac{8\nu_r^2}{\aal} |\Qm u|^2}
\end{array}
\ed
because $(\Pm \rot \om, \Qm u)=0$. In order to estimate the form $b$ we write 
\be\label{dmh40}
b(u,u_1,\Qm u) = b(\Pm u,u_1,\Qm u) + b(\Qm u,u_1,\Qm u)
\ee
and (because $b(u,v,v)=0$)
\be\label{dmh50}
b(u_2,u,\Qm u) = b(u_2,\Pm u,\Qm u).
\ee
Using \eqref{b4} and the Young inequality we infer that
\be \label{dmh60}
\begin{array}{rl}
b(\Pm u,u_1,\Qm u) & \leq c_1 |\Pm u|^{1/2} ||\Pm u ||^{1/2} |u_1|^{1/2} ||u_1||^{1/2} ||\Qm u||, \\ 
b(\Qm u,u_1,\Qm u) &\leq c_1 |\Qm u| \cdot ||\Qm u || \cdot ||u_1|| \\[1.5ex]
				& \leq \ds{\frac{\nu+\nu_r}{8} ||\Qm u ||^2 + \frac{2c_1^2}{\nu+\nu_r} |\Qm u |^2 ||u_1||^2, }\\[1.5ex]
b(u_2,\Pm u,\Qm u) &\leq c_1 |u_2|^{1/2} ||u_2||^{1/2} |\Pm u|^{1/2}||\Pm u ||^{1/2} ||\Qm u||.
\end{array}
\ee

Now we treat the equation \eqref{dmh20} in a similar manner. Taking it's scalar product with  $\QQm \om$ in $\dot H^0 _p$ we obtain
\be \label{dmh70}
\begin{array}{rl}
\ds{\frac{1}{2}\frac{d}{dt}|\QQm \om|^2 + \aal ||\QQm \om ||^2 + b_1(u,\om_1,\QQm \om) }&+ b_1(u_2,\om,\QQm \om) \\
			&= 2\nu_r (\rot u, \QQm \om) + (g,\QQm \om).
\end{array}
\ee
The terms of the RHS of \eqref{dmh70} are estimated as follows
\bd 
\begin{array}{rl}
(g, \QQm \om) &\leq \ds{\frac{\aal}{4}||\QQm \om||^2 + \frac{1}{\aal}||g||_{H^{-1}}^2, }\\[1.5ex]
2\nu_r(\rot u, \QQm \om) & = 2\nu_r (\QQm \rot u, \QQm \om) \leq \ds{ \frac{\nu_r}{4}||\Qm u||^2 + 4\nu_r |\QQm \om|^2    } 
\end{array}
\ed
and the form $b_1$ as will be stated next
\be
\begin{array}{rl}
b_1(\Pm u ,\om_1,\QQm \om) \leq & c_1 |\Pm u|^{1/2}||\Pm u ||^{1/2} |\om_1|^{1/2}||\om_1||^{1/2} ||\QQm\om||, \\
b_1(\Qm u,\om_1,\QQm \om) \leq & c_1|\Qm u|^{1/2}||\Qm u||^{1/2} ||\om_1|| \cdot |\QQm \om|^{1/2} ||\QQm \om||^{1/2} \\[1.5ex]
				\leq &\ds{ \frac{c_1||\om_1||}{2}\lee(|\Qm u|\cdot ||\Qm u|| + |\QQm \om| \cdot ||\QQm \om|| \ri) }\\[1.5ex]
				\leq & \ds{ \frac{\nu+\nu_r}{8} ||\Qm u||^2 + \frac{c_1^2||\om_1||^2}{2(\nu+\nu_r)}|\Qm u|^2 
								+\frac{\aal}{8}||\QQm \om||^2 + \frac{c_1^2||\om_1||^2}{2\aal}|\QQm \om|^2,     }\\
b_1(u_2,\PPm \om,\QQm \om) \leq & c_1 |u_2|^{1/2}||u_2||^{1/2}||\QQm \om|| \cdot |\PPm \om|^{1/2} ||\PPm \om||^{1/2}.
\end{array}
\ee
Adding \eqref{dmh30} and \eqref{dmh70} and using foregoing estimates we arrive at
\be \label{dmh90}
\begin{array}{rl}
\ds{ \frac{d}{dt} \lee(|\Qm u|^2 + |\QQm \om|^2 \ri) }& + k_1 \lee( ||\Qm u||^2 + ||\QQm \om||^2 \ri)  \\
	& \ds{ -\lee(|\Qm u|^2 + |\QQm \om|^2 \ri) \lee(\frac{16\nu_r^2}{\aal} + \frac{4c_1^2}{k_3}\lee(||u_1||^2 + ||\om_1||^2\ri) \ri) \leq \beta(t) }
\end{array}
\ee
where
\bd
\begin{array}{l}
\beta(t) = \textrm{ all terms converging to 0 as } t \to \infty, \\
k_3 = \min(\nu+\nu_r,\aal).
\end{array}
\ed
We make use of the following inequalities $\la_{m+1} |\Qm u|^2 \leq ||\Qm u||^2$ and $\la_{m+1} |\QQm \om||^2 \leq ||\QQm \om||^2$ in order to write \eqref{dmh90} in a form which allows us to use the generalized Gronwall Lemma (Lemma \ref{gron})	
\be \label{dmh100}
\begin{array}{rl}
\ds{ \frac{d}{dt} }\big(|\Qm u|^2 &+ |\QQm \om|^2 \big)   \\
	& +\ds{ \lee( |\Qm u|^2 + |\QQm \om|^2 \ri)	\cdot \lee(k_1\la_{m+1} -\frac{4c_1^2}{k_3} \lee(||u_1||^2 + ||\om_1||^2\ri) - \frac{16 \nu_r^2}{\aal}  \ri)  		\leq \beta}.
\end{array}
\ee
Denoting
\bd
\begin{array}{l}
\xi(t) = |\Qm u|^2 + |\QQm \om|^2,\\[1.5ex]
\ds{\gamma(t) = k_1\la_{m+1} -\frac{4c_1^2}{k_3} \lee(||u_1||^2 + ||\om_1||^2\ri) - \frac{16 \nu_r^2}{\aal}, }
\end{array}
\ed
we can write \eqref{dmh100} in the form
\bd
\frac{d\xi}{dt} + \gamma \xi \leq \beta.
\ed
Now we have only to check the assumptions of Lemma \ref{gron}. In section \ref{ape} we have shown that (see \eqref{est60})
\bd
\frac{1}{T} \int_t ^{t+T} ||u_1(s)||^2 + ||\om_1(s)||^2\, ds \leq \frac{2}{k_1^2} \tf_{-1} ^2.
\ed
To check assumption \eqref{l1} we write
\bd
\begin{array}{rl}
\ds{ \liminf_{t\to \infty} \frac{1}{T} \int _t ^{t+T} \gamma(s)\, ds }  \geq & 
				\ds{k_1 \la_{m+1} - \frac{16 \nu_r^2}{\aal} - \limsup_{t \to \infty} \frac{2c_1^2}{k_3} \lee(||u_1||^2 + ||\om_1||^2 \ri) }\\
			 \geq & \ds{k_1 \la_{m+1} - \frac{16 \nu_r^2}{\aal} - \frac{8c_1^2}{k_1^2 k_3} \tf _{-1}^2 }.
\end{array}
\ed
This assumption is satisfied for
\be \label{dmh110}
m \geq \frac{16 \nu_r^2}{d\la_1\aal k_1}  + \frac{8c_1^2}{d\la_1k_3k_1^3} \tf^{2}_{-1}
\ee
It's easy to check that if $m$ satisfies \eqref{dmh110} then the other assumptions are also fulfilled. That ends the proof.
\end{proof}

\subsection*{Corollaries}
This part of the paper was inspired by the paper of J.C. Robinson \cite{robscales}, in which he shows how the spatial distribution of a force influences the dimension of a global attractor for Navier-Stokes equation.

Suppose that the asymptotic amount of forces $f$  and moments $g$ is equal to $\tf$.
We check how the spatial distribution of external forces and moments influences the number of determining modes.
We consider some cases and write the calculations only for $f$ because these for $g$ are exactly the same.

\begin{enumerate}
\item
Let us consider that the forces and moments act only in two scales i.e. 
\bd
|f_n(t)|^2 = |f_N(t)|^2 =\frac{|f(t)|^2}{2}.
\ed
Then the $\mbb{H}^{-1}$ norm of $f$ satisfies
\bd
||f||_{\mbb{H}^{-1}}^2 =\frac{|f|^2}{2} \lee( \frac{1}{\la_n} + \frac{1}{\la_N} \ri) \sim |f|^2\lee(\frac{1}{n} + \frac{1}{N}\ri).
\ed
That, after inserting to \eqref{dmh110}, gives us 
\bd
m \geq \frac{16 \nu_r^2}{d\la_1\aal k_1}  + \frac{8c_1^2}{d\la_1k_3k_1^3} 
 \tf^2 \lee(\frac{1}{n} + \frac{1}{N}\ri).
\ed
Therefore the number of determining modes depends on inverse of the number of modes in which forces and moments are acting.

\item
Suppose that forces and moments are acting only in some scales and we know nothing about the distribution of it's magnitude through modes
\bd
f= \sum_{k=n}^{N} f_k w_k.
\ed
The following inequalities are straightforward consequence of the definition of the norm in $\mbb{H}^{-1}$ 
\bd
\frac{1}{\la_N} |f|^2 \leq ||f||_{\mbb{H}^{-1}}^2 \leq \frac{1}{\la_n} |f|^2
\ed 
 and we arrive at
\be \label{dmh120}
\frac{1}{\la_N} \tf^2 \leq \tf_{-1}^2 \leq \frac{1}{\la_n} \tf^2.
\ee
Inserting \eqref{dmh120} to \eqref{dmh110} we get
\bd
m \geq \frac{16 \nu_r^2}{d\la_1\aal k_1}  + \frac{8c_1^2}{d\la_1k_3k_1^3}  \tf_{-1}^2
		\geq \frac{16 \nu_r^2}{d\la_1\aal k_1}  + \frac{8c_1^2}{d\la_1k_3k_1^3}  \frac{\tf^2}{\la_N}.
\ed

\item
Assume that the magniude of forces is uniormly distributed in the first $N$ modes that is $|f_k|^2 =(1/N) |f|_{L^2}^2$. Then the following equalities hold 
\bd
\frac{||f||_{\mbb{H}^{-1}}^2}{|f|^2} = \frac{1}{N} \sum_{k=1}^{N} \la_k ^{-1}.
\ed  
Since for large $k$ is $\la_k \sim c\cdot k$ then $\tf_{-1} \sim N^{-1/2}(\ln N)^{1/2} \tf$ and we have the following estimate of the number of determining modes
\bd
m \geq \frac{16 \nu_r^2}{d\la_1k_3k_1}  + \frac{8c_1^2 }{d\la_1k_3k_1^3}
 \tf^2N^{-1/2}(\ln N)^{1/2}.
\ed

\item
Let us consider that the forces and moments act in first $N$  scales and the norm of each mode incerases linearly with it's number 
\bd
|f_k|^2 = \frac{2||f||_{L^2}^2}{N(N+1)} \cdot k.
\ed
for $k=1, \ldots N$. Then we have
\bd
||f||_{H^{-1}}^2 = \sum_{k=1}^{n} \frac{1}{\lambda_k} \frac{2k||f||_{L^2}^2}{n(n+1)} 
		= \frac{2||f||_{L^2}^2}{n(n+1)} \sum_{k=1}^{n} \frac{k}{\lambda_k}.
\ed
Since $\lambda_k \sim k$ then
\bd
||f||_{H^{-1}}^2 \sim \frac{||f||_{L^2}^2}{(n+1)}.
\ed
and the  \eqref{dmh110} implies
\bd
m \geq \frac{16 \nu_r^2}{d\la_1k_3k_1}  + \frac{8c_1^2 }{d\la_1k_3k_1^3}  \frac{c\tf^2}{n+1}.
\ed

\item
Now we assume that the norm of each mode decerases linearly with it's number, that is
\bd
|f_k|^2 = \frac{2||f||_{L^2}^2}{n(n+1)} \cdot (n+1-k).
\ed
Thus we have
\bd
||f||_{H^{-1}}^2 = \frac{2||f||_{L^2}^2}{n(n+1)} \sum_{k=1}^{n} \frac{n+1-k}{\lambda_k} \approx 
 \frac{2||f||_{L^2}^2}{n(n+1)} \sum_{k=1}^{n} \lee( \frac{n+1}{k} -1 \ri)  \approx 
 \frac{2||f||_{L^2}^2}{n(n+1)} [(n+1)\ln n -n].
\ed
From \eqref{dmh110} and the foregoing we infer that
\bd 
m \geq \frac{16 \nu_r^2}{d\la_1 k_3 k_1}  + \frac{8c_1^2 }{d\la_1k_3k_1^3}
 \frac{2\tf^2}{n(n+1)} [(n+1)\ln n -n].
\ed
\end{enumerate}
Above considerations show that if we increase the number of modes in which the forces and moments act or we act only in modes with high-wavenumber then the number of modes necessary to determine the flow decreases. It could be so because in small scales, corresponding to high-wavenumber modes, the damping effect of viscosity is stronger than in large scales. Moreover the number of determining modes depends on how are the forces and moments distributed throughout the modes.

\section{Determining nodes}\label{nodes}

In this section we define the notion of determining nodes and estimate the number of them in terms of asymptotic strength of forces and moments 
$\tf$.

We consider two solutions of micropolar fluid equations $(u_1,\om_1)$ and $(u_2,\om_2)$:
\be \label{rr11}
\frac{\partial u_{1}}{\partial t} - (\nu + \nu_{r}) \D u_1 + (u_1 \cdot \n) u_1 + \n p = 2 \nu_{r} \rot \om_1 + f_1,
\ee
\be\label{rr12}
\di u_1 = 0,
\ee
\be\label{rr13}
\frac{\partial \om _1}{\partial t} - \alpha \D \om_1  + (u_1 \cdot \n) \om_1 +4\nu_r \om_1  = 2 \nu_r \rot u_1 +g_1
\ee
and
\be \label{rr21}
\frac{\partial u_{2}}{\partial t} - (\nu + \nu_{r}) \D u_2 + (u_2 \cdot \n) u_2 + \n q = 2 \nu_{r} \rot \om_2 + f_2,
\ee
\be\label{rr22}
\di u_2 = 0,
\ee
\be\label{rr23}
\frac{\partial \om _2}{\partial t} - \alpha\D \om_2  + (u_2 \cdot \n) \om_2 + 4\nu_r \om_2  = 2 \nu_r \rot u_2+g_2,
\ee
corresponding to two possibly different pairs of external forces and moments $(f_1,g_1)$ and $(f_2,g_2)$ respectively,
the corresponding pressure terms are $p=p(x,t)$ and $q=q(x,t)$. The boundary conditions are periodic with vanishing space averages for both problems. We assume, as in previous section, that forces and moments have the same asymptotic behavior but now we want them to converge in the $L^2$ norm, that is,
\be \label{s2}
\int _Q \lee(|f_1(x,t) - f_2(x,t)|^2 + |g_1(x,t) - g_2(x,t)|^2 \ri)\, dx \to 0  \textrm{ as } t \to \infty.
\ee
We consider a set of $N$ measurement points (called nodes), denoted by $\Sigma$,  $\Sigma = \{x^1,x^2,\ldots,x^N\}$.
We assume, that these points are uniformly distributed within the domain $Q$ in the sense that $Q$ may be covered by $N$ identical squares $Q_1,\ldots,Q_N$ such that there is one and only one given point in each square: $x^i \in Q_i$.

We assume, that the flows have the same time-asymptotic behavior at measurement points. This condition can be written in the form:
\be \label{du} 
\max_{j=1,\ldots,N}|u_1(x^j,t) - u_2(x^j,t)| \to 0, \textrm{ as } t\to \infty 
\ee 
and 
\be \label{dom}
\max_{j=1,\ldots,N}|\om_1(x^j,t) - \om_2(x^j,t)| \to 0, \textrm{as } t\to \infty. 
\ee 
We want to estimate how many points of observation are necessary to determine the asymptotic behavior of a flow in the following sense:

\begin{df}\label{ddn}
The set $\Sigma$ is called a set of \emph{determining nodes} if \eqref{du}, \eqref{dom} together with condition for forces and moments \eqref{s2} implies
\be \label{dn1}
\int _Q \big(|u_1(x,t)-u_2(x,t)|^2+|\om_1(x,t)-\om_2(x,t)|^2\big) \, dx \to 0 \textrm{ as } t \to \infty.
\ee
\end{df}
We will actually show, that \eqref{du}, \eqref{dom} and \eqref{s2} imply that solutions converge to each other in a stronger norm associated with enstrophy, that is
\bd
\int_Q \lee(|\n u_1(x,t) - \n u_2(x,t)|^2 + |\n \om_1(x,t) - \n \om_2(x,t)|^2 \ri)\, dx \to 0 \textrm{ as } t\to\infty.
\ed
Let us denote
\bd
\eta (w) = \max_{1\leq j\leq N} |w(x^j)|
\ed
for each velocity or microrotation field $w$.

There are two lemmas used in the proof of existence of a finite set of determining nodes. One of them is generalized Gronwall's lemma, already used in previous section, the other one is the following lemma from \cite{jones}.
\begin{lm}\label{lm2}
Let the domain $Q$ be covered by $N$ identical squares. Consider the set $\Sigma=\{x^1,\ldots x^N\}$ of points in $Q$, distributed one in each square. Then for each vector field in $\dot{\mbb{H}}_p^2$, the following quantities hold:
\begin{eqnarray}
|w|^2 &\leq& \frac{c}{\lambda_1}\eta(w)^2 + \frac{c}{\lambda_1^2 N^2}|\D w|^2,\\
||w||^2 &\leq& cN\eta(w)^2 + \frac{c}{\lambda_1 N}|\D w|^2,\\
||w||_{L^\infty (Q)}^2&\leq& cN\eta(w)^2 + \frac{c}{\lambda_1 N}|\D w|^2,
\end{eqnarray}
where the constant $c$ depends only on the shape of $Q$ i.e. $L_1/L_2$.
\end{lm}
\begin{tw}
Let $Q$ be domain covered by $N$ identical squares $Q_1,\ldots,Q_N$ and consider a set $\Sigma = \{x^1, \ldots,x^N
\}$ of points in $Q$ distributed one in each square $x^i \in Q_i$, for $1\leq i\leq N$. Let $f_1$ and $f_2$ be two forcing terms in $L^\infty(0,\infty;H)$, $g_1$ and $g_2$ be two moments in $L^\infty(0,\infty; \dot{H}_p ^0)$, satisfying 
\eqref{s2}.

Then the set $\Sigma$ is a set of determining nodes in the sense of Definition \ref{ddn} for the 2-dimensional micropolar fluid equations with periodic boundary conditions provided that
\begin{align*}
N  \geq &\frac{c}{\la_1 k_1} \Bigg\{\frac{8\nu_r^2}{\aal} -2\nu_r + \lee(\frac{c_1^2 c^{1/2}}{\la_1 \nu}
                        + \frac{c_1}{\aal}\ri) \\
            &\cdot \lee(\frac{5\aal k_2 + 32\nu_r^2}{\aal k_1^2 k_2}\tf^2
                + \frac{16C\h{c}_1}{\aal \nu k_1^2 k_2^3} \tf^6 \exp(\h{c}_2 +\h{c}_3 \tf^4) \ri) \\
            & +\frac{16c_1^4\h{c}_1}{\la_1 \aal \nu k_1 k_2}\tf^4 \exp(\h{c}_2 +\h{c}_3 \tf^4 )  \Bigg\},
\end{align*}
where $\tf$ is defined in \eqref{F}.
\end{tw}

\begin{proof}
Let us denote $u=u_1-u_2$, $f=f_1-f_2$ etc. Subtracting equations  \eqref{rr11} and \eqref{rr21} we find 
\be \label{dn10}
\frac{\partial u}{\partial t} + (\nu + \nu_r) \D u + (u_1 \cdot \n) u_1 -(u_2 \cdot \n) u_2 = 2 \nu_r \rot \om + f.
\ee 
By taking the inner product of \eqref{dn10} and $Au$ in H, we get
\be\label{dn20} 
\barray{rl}
\ds{\pdt||u(t)||^2 + (\nu+\nu_r)|Au|^2 }&+b(u,u_1,Au) + b(u_2,u,Au) \\
                                &= 2\nu_r (\rot \om,Au) + (f,Au).
\earray
\ee
Exploiting the orthogonality property \eqref{ortho} we obtain (c.f. \cite{foi})
\bd
b(u,u_1,Au) + b(u_2,u,Au) = -b(u,u,Au_1),
\ed
thus we can write \eqref{dn20} in the form
\be \label{dn30}
\pdt||u(t)||^2 + (\nu+\nu_r)|Au|^2  = 2\nu_r (\rot \om,Au) + (f,Au) + b(u,u,Au_1).
\ee
We estimate the terms in the RHS of \eqref{dn30} using \eqref{b4}, Lemma \ref{lm2} and Young's inequality
\bd
\begin{array}{rl}
2\nu_r(\rot\om,Au) \leq &\ds{\frac{\nu_r}{2}|Au|^2 + 2\nu_r ||\om||^2,  }\\
(f,Au)      \leq & \ds{\frac{\nu_r}{2}|Au|^2 + \frac{1}{2\nu_r}|f|^2,}\\
b(u,u,Au_1)      \leq & c_1|u|^{1/2}||u|| \cdot |Au|^{1/2}|Au_1| \\
            \leq & \ds{\frac{c^{1/4}c_1}{\la_1^{1/4}}\eta(u)^{1/2}||u||\cdot|Au|^{1/2}|Au_1|}\\
                & \ds{+ \frac{c^{1/4}c_1}{\la_1^{1/2}N^{1/2}}|Au|\cdot||u||\cdot|Au_1|   }\\
            \leq & \ds{\frac{c^{1/4}c_1}{\la_1^{1/4}}\eta(u)^{1/2}||u||\cdot|Au|^{1/2}|Au_1| }\\
                &\ds{   + \frac{c^{1/2}c_1^2}{\la_1 N \nu}||u||^2|Au_1|^2  +\frac{\nu}{4}|Au|^2  }.
\end{array}
\ed
Using above estimates we find from \eqref{dn30} that
\be\label{dn40} \barray{rl}
\ds{\pdt ||u||^2 }& \ds{+\frac{3}{4}\nu |Au|^2 - \frac{c^{1/2}c_1^2}{\la_1 N \nu}||u||^2|Au_1|^2  }\\[1.5ex]
        & \ds{\leq    2\nu_r ||\om||^2+ \frac{1}{2\nu_r}|f|^2
            +\frac{c^{1/4}c_1}{\la_1^{1/4}}\eta(u)^{1/2}||u||\cdot|Au|^{1/2}|Au_1|.  }
\earray 
\ee 
Now we treat equations for microrotation \eqref{rr13} and \eqref{rr23} in similar way. Subtracting them we find
\be\label{dn50} 
\frac{\partial \om}{\partial t} + \alpha A_1 \om  + (u_1 \cdot \n)\om_1 - (u_2\cdot \n)\om_2 + 4\nu_r \om  
			= 2 \nu_r \rot u + g. 
\ee 
Multiplying \eqref{dn50} by $A_1\om$ and integrating over $Q$ we get 
\be \label{dn60} 
\barray{rl}
\ds{\pdt ||\om||^2 +\aal|A_1\om|^2 }&+ b_1(u,\om,A_1\om) + b_1(u_2,\om,A_1\om) + 4\nu_r||\om||^2 \\
&=2\nu_r(\rot u,A_1\om) + (g,A_1\om).
\earray
\ee
We estimate the nonlinear terms as follows
\bd
\barray{rl}
b_1(u,\om_1,A_1\om) \leq & c_1 |u|^{1/2}|Au|^{1/2}||\om_1||\cdot |A_1\om| \\
            \leq & \ds{\frac{\aal}{8}|A_1\om|^2 +\frac{2c_1^2}{\aal}|u|\cdot|Au|\cdot||\om_1|^2 }\\[1.5ex]
            \leq &  \ds{\frac{\aal}{8}|A_1\om|^2  + \frac{\nu}{4}|Au|^2 + \frac{4c_1^4}{\aal^2\nu}|u|^2||\om_1||^4,  }\\[1.5ex]
b_1(u_2,\om,A_1\om) \leq & c_1|u_2|^{1/2}|Au_2|^{1/2}||\om||\cdot|A_1\om| \\[1.5ex]
            \leq & \ds{\frac{\aal}{8}|A_1\om|^2 + \frac{2c_1^2}{\aal}|u_2|\cdot|Au_2|\cdot||\om||^2}
\earray
\ed
and the terms in the RHS of \eqref{dn60} similar to the terms in the RHS of \eqref{dn30}, that is
\bd
\barray{rl}
2\nu_r(\rot u,A_1\om) \leq &\ds{\frac{\aal}{8} |A_1\om|^2 + \frac{8\nu_r^2}{\aal}||u||^2, }\\[1.5ex]
(g,A_1\om) \leq & \ds{\frac{\aal}{8} |A_1\om|^2 + \frac{2}{\aal}|g|^2. }
\earray
\ed
Using above estimates we infer from \eqref{dn60} that
\be\label{dn70}
\barray{rl}
\ds{\pdt ||\om||^2 +\frac{\aal}{2}|A_1\om|^2 }& \ds{+ 4\nu_r||\om||^2 -\frac{8\nu_r^2}{\aal}||u||^2
                    -\frac{2c_1^2}{\aal}|u_2|\cdot|Au_2|\cdot||\om||^2  }\\
            & \ds{\leq \frac{2}{\aal}|g|^2 + \frac{\nu}{4}|Au|^2 + \frac{4c_1^4}{\aal^2\nu}|u|^2||\om_1||^4.  }
\earray
\ee
Adding \eqref{dn40} and \eqref{dn70} we have
\be\label{dn75}
\barray{rl}
\ds{\pdt }&\ds{\big(||u||^2+||\om||^2\big) + \frac{k_1}{2}\lee(|Au|^2 + |A_1\om|^2\ri) }\\
    & \ds{  -||u||^2\lee(\frac{c_1^2 c^{1/2}}{\la_1 N \nu}|Au_1|^2+\frac{8\nu_r^2}{\aal}\ri) }\\
    & \ds{  -||\om||^2\lee(\frac{2c_1}{\aal}|u_2|\cdot|Au_2|-2\nu_r \ri) -\frac{4c_1^4}{\aal^2 \nu} |u|^2||\om_1||^4 }\\
    & \ds{\leq \frac{2}{k_3}\lee(|f|^2 + |g|^2\ri) + \frac{c_1c^{1/4}}{\la_1^{1/4}}\eta(u)^{1/2}||u||\cdot|Au|^{1/2}|Au_1|, }
\earray
\ee
where $k_3=\min\{\nu_r,\aal\}$. By Lemma \ref{lm2} we conclude that
\be \label{dn80}
|Au|^2 \geq \frac{\la_1 N}{c}||u||^2 - \la_1 N^2 \eta(u)^2,
\ee
and analogous inequality for $\om$
\be \label{dn90}
|A_1\om|^2 \geq \frac{\la_1 N}{c}||\om||^2 - \la_1 N^2 \eta(\om)^2.
\ee
By the Poincar\'e inequality we get
\be \label{dn100}
\frac{4c_1^4}{\aal^2 \nu}|u|^2 ||\om_1||^4 \leq \frac{4c_1^4}{\aal^2 \nu \la_1}||u||^2 ||\om_1||^4.
\ee
Taking into account \eqref{dn80} -- \eqref{dn100} we infer from \eqref{dn75} that
\be \label{dn110}
\barray{rl}
\ds{\pdt \big(||u||^2 }& \ds{+||\om||^2\big) +  \Bigg( \frac{k_1\la_1 N}{c} -\frac{c_1^2 c^{1/2}}{\la_1 N \nu}|Au_1|^2
								-\frac{8\nu_r^2}{\aal}-\frac{2c_1}{\aal}|u_2|\cdot|Au_2|  }\\
        &\ds{+  2\nu_r -\frac{4c_1^4}{\aal^2 \nu \la_1}||\om_1||^4   \Bigg) \lee( ||u||^2 +||\om||^2 \ri) }\\
    \leq&\ds{ \frac{2}{k_3}\lee(|f|^2 + |g|^2\ri) + \frac{c_1c^{1/4}}{\la_1^{1/4}}\eta(u)^{1/2}||u||\cdot|Au|^{1/2}|Au_1| }\\
        &\ds{+\frac{k_1\la_1 N^2}{2}(\eta(u)^2 + \eta(\om)^2    ).}
\earray
\ee
Denoting
\bd
\barray{rl}
\gamma(t)   =& \ds{2\lee(\frac{k_1\la_1 N}{c} -\frac{c_1^2 c^{1/2}}{\la_1 N \nu}|Au_1|^2-\frac{8\nu_r^2}{\aal}
            -\frac{2c_1}{\aal}|u_2|\cdot|Au_2|+  2\nu_r -\frac{4c_1^4}{\aal^2 \nu \la_1}||\om_1||^4  \ri),}\\
\beta(t)    =& \ds{2\Bigg(\frac{2}{k_3}\lee(|f|^2 + |g|^2\ri)
            + \frac{c_1c^{1/4}}{\la_1^{1/4}}\eta(u)^{1/2}||u||\cdot|Au|^{1/2}|Au_1| +\frac{k_1\la_1 N^2}{2}(\eta(u)^2 + \eta(\om)^2    ) \Bigg), }\\
\xi(t)  =& ||u||^2 + ||\om||^2 
\earray 
\ed 
we can write \eqref{dn110} in the form 
\bd 
\frac{d\xi}{dt} + \gamma \xi \leq \beta. 
\ed 
The time average of $\beta = \beta(t)$ goes to zero as time goes to infinity because for $t$ bounded away from zero, time
averages of the square of norms of $u_1$ and $u_2$ in $D(A)$ are uniformly bounded (see section \ref{ape}), so the condition \eqref{l3} of Lemma \ref{gron} is satisfied. Using again the same argument one can check that assumption 
\eqref{l2} is satisfied.

In order to check assumption \eqref{l1} we write
\begin{align*}
\liminf_{t\to \infty}\frac{1}{T}\int_t^{t+T}\gamma(s)\,ds \geq& 2\Bigg(\frac{k_1 \la_1 N}{c}+2\nu_r-\frac{8\nu_r^2}{\aal}\\
               &-\limsup_{t\to \infty} \frac{1}{T} \int_t^{t+T} \Bigg\{\frac{c_1^2 c^{1/2}}{\la_1 N \nu}|Au_1(s)|^2\\
                &+\frac{2c_1}{\aal}|u_2(s)|\cdot|Au_2(s)|
                + \frac{2c_1^4}{\la_1 \aal \nu}||\om_1(s)||^4 \Bigg\}\, ds \Bigg)\\
            \geq & 2\Bigg(\frac{k_1 \la_1 N}{c} + 2\nu_r - \frac{8\nu_r^2}{\aal} \\
                &-\limsup_{t\to \infty}\frac{1}{T} \int_t^{t+T} \Bigg\{\frac{c_1^2 c^{1/2}}{\la_1 N \nu}|Au_1(s)|^2
                +\frac{c_1}{\aal}|u_2(s)|^2 \\
                &+\frac{c_1}{\aal}|Au_2(s)|^2
                + \frac{4c_1^4}{\la_1 \aal \nu}||\om_1(s)||^4 \Bigg\}\, ds \Bigg) \\
            \geq& 2\Bigg(\frac{k_1 \la_1 N}{c} + 2\nu_r - \frac{8\nu_r^2}{\aal}
             - \lee(\frac{c_1^2 c^{1/2}}{\la_1 \nu} + \frac{c_1}{\aal}\ri)\\
                &\cdot \lee( \frac{5\aal k_2 + 32\nu_r^2}{\aal k_1^2 k_2}\tf^2
                + \frac{16C\h{c}_1}{\aal \nu k_1^2 k_2^3} \tf^6 \exp(\h{c}_2 +\h{c}_3 \tf^4) \ri) \\
            & -\frac{16c_1^4\h{c}_1}{\la_1 \aal \nu k_1 k_2}\tf^4 \exp(\h{c}_2 +\h{c}_3 \tf^4 ) \Bigg).
\end{align*}
Therefore if
\be
\begin{array}{rl}\label{dn120}
N  \geq &\ds{\frac{c}{\la_1 k_1} \Bigg\{ \frac{8\nu_r^2}{\aal} -2\nu_r + \lee(\frac{c_1^2 c^{1/2}}{\la_1 \nu}
                        + \frac{c_1}{\aal}\ri)  }\\
            &\ds{\cdot \lee(\frac{5\aal k_2 + 32\nu_r^2}{\aal k_1^2 k_2}\tf^2
                + \frac{16C\h{c}_1}{\aal \nu k_1^2 k_2^3} \tf^6 \exp(\h{c}_2 +\h{c}_3 \tf^4) \ri) }\\
            &\ds{ +\frac{16c_1^4\h{c}_1}{\la_1 \aal \nu k_1 k_2}\tf^4 \exp(\h{c}_2 +\h{c}_3 \tf^4 )  \Bigg\}  },
\end{array}
\ee
then \eqref{l1} holds, other assumptions are also fulfilled and we infer from uniform Gronwall lemma that $||u(t)||^2 + ||\om(t)||^2 \to 0$ as $t \to \infty$. That ends the proof.
\end{proof}

\section{Hausdorff and fractal dimension of attractor}\label{dimension}
In this section we recall the notions of Hausdorff and fractal dimension of an invariant set. Then we show that  the attractors $\mc{A}_{\nu_r}, \nu_r \geq 0$ associated with the system \eqref{eq1}-\eqref{eq3} have finite Hausdorff and fractal dimensions, which can be estimated by constants depending on the data: $f,g,\nu,\alpha$ and the domain of flow $Q$ but independent of the microrotation viscosity $\nu_r$.

Let $(X,d)$ be a metric space and $Y \subset X$ be a subset of $X$. For every $d \in \mathbb{R}_+$ and $\epsilon >0$ we set
\begin{displaymath}
\mu_H (Y,d,\epsilon) = \inf \sum _{i \in I} r_i ^d,
\end{displaymath}
where the infimum is for all coverings of $Y$ by a family $(B_i)_{i \in I}$ of balls of radii $r_i \leq \epsilon$. $\mu_H (Y,d,\epsilon)$ is non increasing function of $\epsilon$.
The number $\mu_H (Y,d)$, called the $d$-dimensional Hausdorff measure of set $Y$, is defined as
\begin{displaymath}
\mu_H (Y,d)= \lim _{\epsilon \to 0} \mu_H (Y,d,\epsilon) = \sup _{\epsilon >0} \mu_H (Y,d,\epsilon).
\end{displaymath}
If $\mu_H (Y,d') < \infty$ for some $d'$ then $\mu_H (Y,d) =0$ for every $d>d'$. Then there exist $d_0 \in [0,\infty]$ such that $\mu_H (Y,d)=0$ for every $d>d_0$ and $\mu_H (Y,d)=\infty$ for $d<d_0$. The number $d_0$ is called the Hausdorff dimension of $Y$ -- $d_H(Y)$.

Now we define the fractal dimension of  $Y$. Let $\epsilon >0$, we denote by $n_Y (\epsilon)$ the minimum number of balls of $X$ of radius $\epsilon$ necessary to cover the set $Y$. The fractal dimension of $Y$, called also the capacity of  $Y$, is the number $d_F$ defined as follows
\begin{displaymath}
d_F (Y) = \limsup _{\epsilon \to 0} \frac{\log n_Y (\epsilon)}{\log 1/\epsilon}.
\end{displaymath}
We refer the reader to \cite{fal} for more details.
\begin{tw} \label{AttrDim}
There exist a constant $C_0$ depending only on the shape of $Q$, which has following property. If $N$ is an integer such that 
\be\label{N}
N-1 < 2C_0 (k_1^3 k_2) ^{-1/2} \left(|f|^2+|g|^2\right)^{1/2} \leq N, 
\ee 
where $k_1,k_2$ are as in \eqref{k}, then the $N$-dimensional volume element in phase space $\mc{H}$ is exponentially decaying, moreover the Hausdorff dimensions of $\mc{A}_{\nu_r}$, $\nu_r \geq 0$ are less or equal to $N$ and the fractal dimensions are less or equal to $2N$.
\end{tw}

\begin{proof}
The relation \eqref{N} proceeds from estimates for the uniform Lyapunov exponents associated with attractors
$\mc{A}_{\nu_r}$.

First, we rewrite our equations \eqref{eq1} - \eqref{eq3} in a more suitable form. In order to do this we introduce some notations.
Let $\bu_i = (u_i,\om_i) \in \mc{H}$ (or $V$) for $i=1,2$. We introduce scalar products and norms in $\mc{H}$ and $\mc{V}$ as follows
\bd
[\bu_1, \bu_2] = (u_1,u_2) + (\om_1,\om_2),\quad [\bu] =[\bu,\bu]^{1/2}
\ed
for all $\bu,\bu_1,\bu_2 \in \mc{H}$ and
\bd
[[\bu_1, \bu_2]] = (\n u_1,\n u_2) + (\n\om_1,\n\om_2), \quad [[\bar{u}]]=[[\bar{u},\bar{u}]]^{1/2}
\ed
for all $\bu,\bu_1,\bu_2 \in \mc{V}$. The notation seems to be confusing, but it will be always clear from context if
$(\cdot,\cdot)$ denotes the scalar product in $L^2$ or a vector in $\mc{H}$ or $\mc{V}$.

We define the trilinear form $B$ on $\mc{V} \times \mc{V} \times \mc{V}$ by
\bd
B(\bu_1,\bu_2,\bu_3) = b(u_1,u_2,u_3) + b_1(u_1,\om_2,\om_3)
\ed
and associate with $B$ the bilinear continuous operator $\mc{B}$ from $\mc{V} \times \mc{V}$ to $\mc{V}'$ as follows
\bd
\langle \mc{B}(\bu_1,\bu_2),\phi\rangle = B(\bu_1,\bu_2,\phi) \quad \bu_1,\bu_2,\phi \in \mc{V} .
\ed

We define bilinear forms $R$ and $a$ from $\mc{V} \times \mc{V}$ by
\bd
\begin{array}{l}
R(\bu_1,\bu_2) = -2\nu_r(\rot \om_1,u_2) - 2\nu_r(\rot u_1,\om_2) + 4\nu_r(\om_1,\om_2), \\
a(\bu_1,\bu_2) = (\nu+\nu_r)(\n u_1,\n u_2) + \alpha (\n \om_1,\n\om_2)
\end{array}
\ed
and associate linear, continuous operators on $\mc{V}$ to $\mc{V}'$ by
\bd
\langle \mc{R}(\bu),\phi\rangle= R(\bu,\phi),
\quad \langle \mc{A}(\bu),\phi\rangle = a(\bu,\phi),
\quad \bu,\phi \in \mc{V}.
\ed
Let $G=(f,g) \in \mc{H}$, then the weak form of \eqref{eq1}-\eqref{eq3}, which has form (c.f. \cite{szopae})
\be \label{dr1}
\begin{array}{rl}
\ds{\frac{d}{dt}(u(t),\varphi(t)) + (\nu + \nu_{r})(\n u(t),\n \varphi(t))} & \ds{+ b(u(t),u(t),\varphi)} \\[1.5ex]
        & \ds{= 2 \nu_{r}(\rot \om(t),\varphi) + (f,\varphi)}
\end{array}
\ee
for all $\varphi \in V$,and
\be \label{dr2}
\begin{array}{rl}
\ds{\frac{d}{dt}(\om(t),\psi) + \alpha(\n \om(t),\n \psi)}& \ds{+b_1 (u(t),\om(t),\psi) +4\nu_r(\om(t),\psi)} \\[1.5ex]
        &\ds{ = 2 \nu_r(\rot u,\psi) +(g(t),\psi)}
\end{array}
\ee
for all $\psi \in \dot H^1 _p (Q)$, in the sense of scalar distributions on $(0,\infty)$, can be rewritten as
\bd
\frac{d}{dt}[\bu,\phi] + a(\bu,\phi) + B(\bu,\bu,\phi) + R(\bu,\phi) = [G,\phi],
\ed
where $\bu=(u,\om), \phi=(\varphi,\psi)$ or in the functional form
\bd
\frac{d}{dt} \bu + \mc{A}(\bu) + \mc{B}(\bu,\bu) +\mc{R}(\bu) = G.
\ed
We will concern ourself with the corresponding, linearized about $\bu$ problem
\bd
\frac{d}{dt}U=F'(\bu)U,
\ed
where $F'(\bu)=-\mc{A}(U) - \mc{B}(\bu,U) - \mc{B}(U,\bu) +\mc{R}(U)$.

Our aim now is to estimate form above the trace
\bd
\tr F'(\bu)\circ \Pn = \sum_{i=1}^{N}[F'(\bu)\ph_j,\ph_j],
\ed
where $\ph_j=\ph_j(\tau)$, $j=1,\ldots,N$ is an orthonormal in $\mc{H}$ basis of $\Pn(\tau)\mc{H}=\textrm{Span}\{U_1(\tau), \ldots,U_N(\tau) \} $, $\Pn(\tau,\xi_1,\ldots,\xi_N)$ being the orthogonal projector in $\mc{H}$ on the space spanned by $U_j(\tau)$, $j=1,\ldots,N$, where $U_j$ satisfy
\bd
\frac{d}{dt}U_j=F'(\bu)U_j, \quad U_j(0) = \xi_j, \quad j=1,\ldots,N.
\ed
Because $B(\bu,\ph_j,\ph_j)=0$, we have
\bd
[F'(\bu)\ph_j,\ph_j]=-a(\ph_j,\ph_j)-B(\ph_j,\bu,\ph_j)-R(\ph_j,\ph_j).
\ed
Let us write $\bu=(u,\om)$ and $\ph_j=(v_j,z_j)$, then we obtain
\begin{align*}
-\sum_{i=1}^{N}&a(\ph_j,\ph_j) = -(\nu+\nu_r)\sum_{i=1}^{N}||v_j||^2 - \alpha \sum_{i=1}^{N}||z_j||^2,\\
-\sum_{i=1}^{N}&B(\ph_j,\bu,\ph_j) = -\sum_{i=1}^{N}b(v_j,u,v_j) -\sum_{i=1}^{N}b_1(v_j,\om,z_j)\\
\intertext{and}
-\sum_{i=1}^{N}&R(\ph_j,\ph_j) = 2\nu_r\sum_{i=1}^{N}(\rot z_j,v_j) +2\nu_r\sum_{i=1}^{N}(\rot v_j,z_j)
        -4\nu_r \sum_{i=1}^{N}|z_j|^2.
\end{align*}
Now let us consider the operator $a+R$ \bd
\begin{array}{rl}
a(\ph_j,\ph_j) + R(\ph_j,\ph_j) =& (\nu+\nu_r)((v_j,v_j))^2 + \alpha((z_j,z_j))^2 \\
            &-2\nu_r(\rot z_j,v_j) -2\nu_r(\rot v_j,z_j) +4\nu_r |z_j|^2.
\end{array}
\ed
Because of \eqref{id1}, Schwartz and Young's inequalities we have
\bd
2\nu_r(\rot z_j,v_j) +2\nu_r(\rot v_j,z_j) = 4\nu_r(\rot z_j,v_j) \leq \nu_r ||v_j||^2 + 4\nu_r|z_j|^2,
\ed
therefore
\bd
a(\ph_j,\ph_j) + R(\ph_j,\ph_j) \geq \nu ((v_j,v_j)) + \alpha ((z_j,z_j)) \geq k_1[[\ph_j,\ph_j]] \quad \forall \ph_j \in V,
\ed
where $k_1=\min\{\nu,\alpha\}$ as in \eqref{k}, whence
\bd
-\sum_{i=1}^{N} \big(a(\ph_j,\ph_j) + R(\ph_j,\ph_j) \big) \leq -k_1 \sum_{i=1}^{N} [[\ph_j,\ph_j]]^2.
\ed
We estimate the trilinear form $b$ as follows
\bd
\left|\sum_{i=1}^{N}b(v_j,u,v_j) \right| = \left| \sum_{i=1}^{N} \int_Q (v_j \cdot \n)u v_j \,dx \right|
            \leq \int_Q |\n u(x,t)| \rho_1(x,t) \,dx,
\ed
where
\bd
\rho_1(x,t) = \sum_{i=1}^{N} |v_j(x,t)|^2,
\ed
and similarly the form $b_1$
\bd
\begin{array}{rl}
\ds{\left|\sum_{i=1}^{N}b_1(v_j,u,z_j) \right| }&=\ds{ \left| \sum_{i=1}^{N} \int_Q (v_j \cdot \n)u z_j \,dx \right| }\\
            &\leq \ds{\int_Q |\n u(x,t)| \rho_1(x,t)^{1/2} \rho_2(x,t)^{1/2} \,dx,}
\end{array}
\ed
where
\bd
\rho_1(x,t) = \sum_{i=1}^{N} |z_j(x,t)|^2.
\ed
Therefore we can estimate the form $B$ as follows
\bd
\left|\sum_{i=1}^{N} B(\ph_j,\bu,\ph_j) \right|
                        \leq \int_Q \left( |\n u|\rho_1 + |\n\om|\rho_1^{1/2} \rho_2^{1/2} \right)\, dx.
\ed
Setting $\rho = \rho_1 + \rho_2$ and by Cauchy and Schwartz inequalities we get
\begin{align*}
\left|\sum_{i=1}^{N} B(\ph_j,\bu,\ph_j) \right| &\leq \int_Q\left( \rho |\n u|+ \rho |\n \om|\right)\,dx \\
    & \leq \left( \int_Q \rho^2\, dx \right)^{1/2} \left(\int_Q (|\n u| + |\n \om|)^2\, dx \right)^{1/2}\\
    & \leq |\rho|  \left(2\int_Q\lee( |\n u|^2 + |\n \om|^2 \ri)\, dx \right)^{1/2}\\
    & = \sqrt{2} |\rho| \cdot [[\bu]].
\end{align*}
From the above estimates we infer that
\be \label{atr5}
\tr F'(\bu(\tau))\circ \Pn(\tau) \leq -k_1 \sum_{i=1}^{N}[[\ph_j(\tau)]]^2 + \sqrt{2} |\rho(\tau)|\cdot [[\bu(\tau)]].
\ee
Because the family $\{\ph_j(\tau)\}_{j=1}^{N}$ is orthonormal in $\mc{H}$, then the family of corresponding pairs
$(v_j,z_j)$ is orthonormal in $L^2(Q)^{2 \times 2}$ and we can use a generalization of the Sobolev-Lieb-Thirring inequality (\cite{tem2}) and write
\be \label{atr10}
|\rho(\tau)|^2 \leq C_0 \sum_{i=1}^{N} \left(||v_j||^2 + ||z_j||^2\right) = C_0 \sum_{i=1}^{N} [[\ph_j(\tau)]]^2,
\ee
where the constant $C_0$ depends only on the shape of $Q$.

Let us perceive that $\int_Q \rho(x,t)\, dx=N$, then by Schwartz inequality we get \bd N^2 \leq |Q|\cdot |\rho|^2, \ed and taking into account \eqref{atr10} 
\be\label{atr20}
\sum_{i=1}^{N}[[\ph_j(\tau)]]^2 \geq \frac{N^2}{C_0|Q|}. 
\ee 
We estimate the second term of the RHS of \eqref{atr5} using Young's inequality, \eqref{atr10} and \eqref{atr20} 
\bd 
\sqrt{2} |\rho|\cdot [[\bu]] \leq \frac{k_1}{2C_0}|\rho|^2 +\frac{C_0}{k_1}[[\bu]]^2
        \leq \frac{k_1}{2}\sum_{i=1}^{N}[[\ph_j]]^2 +  \frac{C_0}{k_1}[[\bu]]^2.
\ed
From the above considerations, we have
\be \label{atr30}
\tr F'(\bu(\tau))\circ \Pn(\tau) \leq -\frac{k_1}{2} \frac{N^2}{C_0|Q|} + \frac{C_0}{k_1}[[\bu]]^2.
\ee
Let $\bu_0=(u_0,\om_0) \in \mc{A}$, $\bu(\tau) = S(\tau)\bu_0$. Let us set
\bd
\begin{array}{c}
\ds{q_N(t) = \sup_{\bu_0 \in \mc{A}} \sup \left\{ \frac{1}{t} \int_0^t \tr F'(S(\tau)\bu_0) \circ \Pn(\tau)\, d\tau
                 \colon \xi_i \in H,\,[\xi_i] \leq 1,\,i=1,\ldots,N \right\}, }\\
\ds{ q_N=\limsup_{t\to \infty} q_N(t)}
\end{array}
\ed
In view of \eqref{atr30}, we conclude that
\be\label{atr40}
q_N \leq -\frac{k_1}{2C_0|Q|}N^2+\frac{C_0}{k_1} \gamma,
\ee
where
\bd
\gamma = \lim_{t\to\infty} \sup_{\bu\in \mc{A}} \frac{1}{t}\int_0^t[[S(\tau)\bu_0]]^2 \,d\tau.
\ed
We can estimate $\gamma$ in terms of the data using \eqref{e5}. 
Since
\bd 
\frac{1}{t} \int_0^t [[S(t)\bu_0]]^2 \leq \frac{1}{k_1 k_2} [G]^2 + \frac{1}{t} \frac{1}{k_1}[\bu_0]^2, 
\ed
then 
\bd 
\gamma \leq \frac{1}{k_1 k_2} [G]^2 
\ed 
and by \eqref{atr40} 
\be\label{atr50} 
q_N \leq -\frac{k_1}{2C_0|Q|}N^2+\frac{C_0}{k_1^2 k_2} [G]^2. 
\ee 
Setting
\bd 
\kappa_1=\frac{k_1}{2C_0|Q|},\qquad \kappa_2=\frac{C_0}{k_1^2 k_2} [G]^2, 
\ed 
we can write \eqref{atr50} in the form 
\bd 
q_N \leq -\kappa_1 N^2 + \kappa_2. 
\ed 
The general theory provided in \cite{tem2} allows us to conclude that the uniform Lyapunov exponents $\mu_j$ associated with the attractor satisfy 
\bd
\mu_1+\cdots+\mu_j \leq -\kappa_1 j^2 + \kappa_2, \qquad \forall j\in \mbb{N}, 
\ed 
and for $N$ satisfying (c.f.. Lemma VI, 2.2 in \cite{tem2}) 
\bd 
N-1 < \left(\frac{2\kappa_2}{\kappa_1}\right)^{1/2} \leq N, 
\ed 
the thesis of Theorem \ref{AttrDim} holds.
\end{proof}

The estimate of the dimension of the global attractor, which we obtained, is similar to analogous estimate for micropolar fluid equations with no-slip boundary conditions (c.f. \cite{luk1}), but of the higher order than analogous estimate for NSE with periodic boundary conditions.

\section{Conclusions}
The lack of orthogonality property of the form $b_1$, which was mentioned is subsection \ref{spaces}, causes that the estimates of numbers of determining modes, nodes and the dimension of global attractor are more involved than corresponding estimates for Navier-Stokes equation in space-periodic case. The reason is following: estimates of the square of the norm of solution in $\mc{V}$ and the average of the square of the norm in $D(A)$ for NSE, which are analogous to \eqref{e7a} and \eqref{Au-norm}, are square with respect to $F$ (c.f. \cite{foi}), where  
\bd 
F =\limsup_{t\to \infty} \lee(\int_Q |f(x,t)|^2 \,dx \ri)^{1/2}.
\ed
The estimates \eqref{e7a} and \eqref{Au-norm} are exponential with respect to $\tf$, which implies exponential estimate of the number of determining nodes. 
We can't obtain linear dependence of the number of determining modes on $\tf$ because the form $b_1$ doesn't have the orthogonality property \eqref{ortho}.

We check how the estimates we obtained depend on $\nu_r$. The case when $\nu_r$ is small  is particularly interesting because if $\nu_r=0$ the micropolar fluid system reduces to the Navier-Stokes system and the velocity field $u$ becomes independent on microrotation field $\om$
\bd
\begin{array}{c}
u_t-\nu \D u+(u\cdot \n)u + \n p =f,  \\
\di u=0, \\
\om_t -\aal \D \om + (u\cdot \n)\om =g.
\end{array}
\ed

\begin{enumerate}
\item The estimate of the number of determining modes $m \approx c_1 \tf_{-1} ^2 + c_2$ (c.f. \eqref{dmh110}) is similar to that for Navier-Stokes equation (\cite{foi}) and micropolar fluid equations (\cite{szopad}), both with Dirichlet boundary conditions. The coefficient $c_1$ doesn't depend on $\nu_r$ but if $\nu_r \to 0$ then $c_2 \to 0$ so, if $\nu_r=0$  our estimate agrees with that for Navier-Stokes equation.

\item The estimate of the number of determining nodes $m \approx P(\tf) + Q(\tf)\exp(R(\tf))$ (c.f. \eqref{dn120}), where $P,Q$ and $R$ are polynomials is similar to that for NSE with Dirichlet boundary conditions (\cite{foi}). If $\nu_r \to 0$ the above estimate remains exponential with respect to $\tf$. Our result is much worse than analogous result for Navier-Stokes equation in \cite{foi}, where it was shown that the dependence is linear.

\item The estimate of the dimension of a global attractor \eqref{N} doesn't depend on $\nu_r$ and is similar to such estimates for NSE \cite{foi} with Dirichlet boundary conditions and for micropolar fluid equations with Dirichlet boundary condition \cite{luk1}
\end{enumerate}



\end{document}